\documentclass[11pt]{article}

\usepackage[T1]{fontenc}
\usepackage[utf8]{inputenc}
\usepackage{mathptmx}

\usepackage{amsmath,amssymb,amsthm,mathtools}
\usepackage{booktabs}
\usepackage{enumitem}
\usepackage{geometry}
\usepackage{microtype}
\usepackage{xcolor}
\usepackage{authblk}
\usepackage{orcidlink}
\usepackage{mathrsfs}
\usepackage{hyperref}
\usepackage[nameinlink,capitalise]{cleveref}

\hypersetup{
  colorlinks=true,
  linkcolor=blue!45!black,
  citecolor=blue!45!black,
  urlcolor=blue!45!black,
  pdftitle={Quantitative Haj\'os Theory and Sharp Stability of Lattice Cube Tilings},
  pdfauthor={Domenico Frijio}
}

\newtheorem{theorem}{Theorem}[section]
\newtheorem{proposition}[theorem]{Proposition}
\newtheorem{lemma}[theorem]{Lemma}
\newtheorem{corollary}[theorem]{Corollary}
\theoremstyle{definition}

\newtheorem{remark}[theorem]{Remark}

\numberwithin{equation}{section}

\DeclareMathOperator{\dist}{dist}
\DeclareMathOperator{\covol}{covol}
\DeclareMathOperator{\LinExt}{LinExt}
\DeclareMathOperator{\lct}{lct}
\DeclareMathOperator{\mult}{mult}
\DeclareMathOperator{\link}{link}
\DeclareMathOperator{\ord}{ord}
\DeclareMathOperator{\sinc}{sinc}
\DeclareMathOperator{\vol}{vol}
\newcommand{\R}{\mathbb R}
\newcommand{\Z}{\mathbb Z}
\newcommand{\C}{\mathbb C}
\newcommand{\one}{\mathbf 1}
\newcommand{\cT}{\mathcal T}
\newcommand{\cH}{\mathcal H}
\newcommand{\cE}{\mathcal E}
\newcommand{\cY}{\mathcal Y}
\newcommand{\cA}{\mathfrak A}

\begin{document}

\title{\fontsize{16}{18}\selectfont\textbf{Quantitative Haj\'os Theory and Sharp Stability of Lattice Cube Tilings}}

\author{Domenico Frijio\orcidlink{0009-0005-2747-3961}}

\affil{Horizon Research, Munich, Germany\\
{\fontsize{9}{10}\selectfont\texttt{dfrijio@horizonrsc.com}}}

\date{\normalfont\fontsize{10}{12}\selectfont 25 August 2026}

\maketitle

\begin{abstract}
Let \(Q_n=[-1/2,1/2]^n\), let \(X_n\) be the moduli space of unimodular lattices in \(\mathbb R^n\), and define cube-overlap energy by summing the covariogram of \(Q_n\) over nonzero lattice vectors. This continuous proper exhaustion vanishes exactly on the compact lattice cube-tiling locus \(\mathcal T_n\). The principal result is the sharp quantitative Haj\'os theorem
\(\operatorname{dist}_{X_n}(\Lambda,\mathcal T_n)\le C_n\mathcal E_n(\Lambda)^{1/n}\)
on a small energy sublevel in every dimension \(n\ge2\), with optimal exponent \(1/n\). The proof uses a fixed-window multiplicity discrepancy. After an endpoint-order subdivision it is a positive sum of rectangular atom volumes, each a product of at most \(n\) affine gaps. Hypergraph blocker duality and Hoffman's error bound yield the required \(n\)-th power estimate.
The paper also develops the global and local structure of \(\mathcal T_n\), a union of \(n!\) compact unipotent nilmanifolds. 
At every tiler, covolume-one coordinate sections form the ideal lattice of a canonical incidence poset. The incident Haj\'os branches are exactly the linear extensions of this poset. Every finite poset occurs in this way.
A squarefree rooted-circuit ideal governs the reduced branch scheme: its multiplicity counts linear extensions, while its log-canonical threshold equals its codimension. At every tiler, overlap is two-sided comparable with the rooted-circuit gauge. A block-squarefree rigidity argument eliminates signed contact cancellations and identifies the optimal local exponent as the reciprocal of the longest rooted circuit. At the standard lattice this sharpens to an exact directed-cycle normal form. Finally, Parseval identifies overlap with a dual sinc-square energy, whose normalized mass escapes every finite window along sharp cyclic families.
\end{abstract}



\section{Introduction}

The theorem of Minkowski and Haj\'os classifies unimodular lattices whose translates of the unit cube tile Euclidean space. After a permutation of the coordinate axes and an integral change of lattice basis, every lattice cube tiler admits a unitriangular basis. The present problem is quantitative: how rapidly must a lattice approach the complete Haj\'os locus when the total overlap of its cube translates tends to zero?

The basic observable is the deck-overlap energy
\begin{equation}\label{eq:intro-energy}
 \cE_n(\Lambda)
 =\sum_{v\in\Lambda\setminus\{0\}}
   \prod_{j=1}^n(1-|v_j|)_+.
\end{equation}
It is intrinsic to the lattice, nonnegative, and vanishes precisely at cube tilers. It simultaneously records the variance of the cube multiplicity on the lattice torus, the complete local contact arrangement, and a dual sinc-square spectrum.

At its core, this work does more than make Haj\'os's theorem quantitative:
it identifies the contact geometry that governs quantitative rigidity.
The overlap energy has the lattice cube-tiling locus as its zero set.
This locus decomposes into finitely many compact unitriangular branches.
Their local incidence is encoded by a canonical poset. The rate at which an
almost-tiling lattice approaches this locus is controlled by the way these
branches meet.

At a tiler $\Lambda$, rooted circuits of $P_\Lambda$ encode the transverse
contacts. The optimal local exponent is exactly $1/c(P_\Lambda)$, where
$c(P_\Lambda)$ is the maximum rooted-circuit length. Since
$c(P_\Lambda)\le n$ at every tiler and a Hamiltonian cycle at the standard
lattice realizes $c=n$, the sharp global exponent is $1/n$. Thus the global
stability law and the stratum-dependent local laws are two scales of the
same contact geometry.

The use of deck overlap as a geometric observable was inspired by the geometric viewpoint introduced in \emph{Global Embedding Limits for Extended Domains in Compact Riemannian Quotients} of D. Frijio~\cite{Frijio}.

The exact theory begins with Haj\'os's proof of Minkowski's conjecture \cite{Hajos}; modern formulations and related cube-tiling results appear in \cite{Kolountzakis1998,Kolountzakis2004,KisielewiczPrzeslawski,LeePfanderWalnut,SteinSzabo}. For an application of Haj\'os's theorem to measure estimates in the geometry of numbers, see Kleinbock, Str\"ombergsson, and Yu \cite{KleinbockStrombergssonYu}.
Fourier characterizations of tiling and spectrality are developed in \cite{Fuglede,IosevichPedersen,Kolountzakis2004}. Stability questions in which the domain, rather than the lattice, varies near a cube belong to a complementary theory \cite{KolountzakisLaba}. 
The exact classification is used here as an input. 
Within the quantitative theory developed in this paper, the following are new: (i) the sharp stability law; (ii) the fixed-window atom method; (iii) the incidence-poset geometry; (iv) the rooted branch scheme; (v) the block-squarefree poset normal form; and (vi) the spectral escape mechanism.

Set
\[
 X_n=\mathrm{SL}_n(\R)/\mathrm{SL}_n(\Z),
\]
equip it with any smooth Riemannian metric, and denote by \(\cT_n\) the locus of lattice cube tilers. The central theorem is the following.

\begin{theorem}[Sharp quantitative Haj\'os theorem]\label{thm:main}
For every \(n\ge2\), there are constants \(C_n,\delta_n>0\) such that
\begin{equation}\label{eq:intro-main}
 \dist_{X_n}(\Lambda,\cT_n)
 \le C_n\cE_n(\Lambda)^{1/n}
\end{equation}
whenever \(\cE_n(\Lambda)\le\delta_n\). No exponent larger than \(1/n\) is valid on every sufficiently small energy sublevel.
\end{theorem}

The new point is the all-strata upper bound. Local algebra at the standard lattice already exhibits Hamiltonian-cycle arcs with distance of order \(|t|\) and energy of order \(|t|^n\), so sharpness is forced. At a nonstandard tiler, however, affine boundary-row forms may cancel before their positive parts are taken; a termwise comparison with rooted circuit monomials is therefore unavailable. The proof of \cref{thm:main} avoids this obstruction by returning to physical space. On a fixed axis-parallel window, the integer multiplicity function is partitioned into rectangular atoms. Each atom volume is a product of at most \(n\) nonnegative affine endpoint gaps. A positive-product error bound, proved from blocker duality and Hoffman's theorem \cite{Hoffman}, controls the \(n\)-th power of the distance to the common zero set. Periodicity and the integer-valued variance identity compare this window discrepancy to \(\cE_n\). Compactness then globalizes the local bounds.

Several structural results accompany the stability theorem. First, \(\cE_n\) is continuous and proper. The explicit cusp estimate
\begin{equation}\label{eq:intro-cusp}
 \cE_n(\Lambda)\ge 2^{1-n}
 \left\lfloor\frac{1}{2\lambda_\infty(\Lambda)}\right\rfloor,
 \qquad
 \lambda_\infty(\Lambda)=\min_{0\ne v\in\Lambda}\|v\|_\infty.
\end{equation}
prevents almost tilers from escaping into the cusp by Mahler compactness. The exact tiling locus is a finite union of compact nilmanifold branches
\begin{equation}\label{eq:intro-branches}
 \cT_n=\bigcup_{\sigma\in S_n}\cY_\sigma,
 \qquad
 \cY_\sigma\simeq N_n(\R)/N_n(\Z),
\end{equation}
where \(N_n\) denotes the upper unitriangular group.

Second, branch incidence at an arbitrary tiler is controlled by a canonical poset. For a coordinate set \(S\subset[n]\), let \(E_S\) be the corresponding coordinate subspace. A set is tight when \(\Lambda\cap E_S\) has full rank and covolume one. Tight sets form the ideal lattice of a unique poset \(P_\Lambda\), and the incident unitriangular flags are exactly its linear extensions. Every finite poset occurs for an explicit algebraic cube-tiling lattice.

Third, the reduced local branch scheme is encoded by a squarefree rooted-circuit ideal \(I_P\). Its facets are the forward tournaments attached to the linear extensions of \(P\). Consequently,
\begin{equation}\label{eq:intro-invariants}
 \mult_0(I_P)=|\LinExt(P)|,
 \qquad
 \lct_0(I_P)=\binom n2.
\end{equation}
After the \(n-1\) independent diagonal equations are adjoined on the determinant-one slice, the log-canonical threshold equals the codimension. Thus multiplicity can encode the \(\#P\)-complete linear-extension count \cite{BrightwellWinkler}, while the birational threshold is universal. The Stanley--Reisner and multiplier-ideal ingredients are classical \cite{Howald,MillerSturmfels}; their realization by cube-tiling branch singularities is the geometric content.

At the standard lattice the overlap itself admits a sharper normal form. For \(B=(b_{ij})\in\mathrm{SL}_n(\R)\) near \(I\), set
\begin{equation}\label{eq:intro-gauge}
 \Delta_n(B)=\sum_{0\ne k\in\Z^n}
 \prod_{i=1}^n(1-|(Bk)_i|)_+,
 \qquad
 \cA(B)=\sum_i|b_{ii}-1|+\sum_C\left|\prod_{(i,j)\in C}b_{ij}\right|,
\end{equation}
where \(C\) ranges over directed simple cycles. Then
\begin{equation}\label{eq:intro-normal}
 \Delta_n(B)=(1+o(1))\cA(B).
\end{equation}
The constant one is uniform across all relative rates and sign chambers. Exact bottleneck duality gives the local exponent \(1/n\); on a prescribed shear digraph \(G\), the optimal exponent is \(1/c(G)\), where \(c(G)\) is the directed circumference. At a general tiler the exact optimal local exponent is first expressed as a rational signed contact-valuation ratio. The fixed-window theorem bounds every such ratio by \(n\), and a block-squarefree refinement proves the poset normal form \(\Delta_U\asymp G_{P_U}\) and the identity \(\kappa(U)=c(P_U)\). Thus the optimal local exponent is \(1/c(P_U)\), while the Hamiltonian cycle shows that the global worst case remains \(1/n\).

Finally, Parseval periodization gives
\begin{equation}\label{eq:intro-spectral}
 \cE_n(\Lambda)=
 \sum_{0\ne\xi\in\Lambda^*}
 \prod_{j=1}^n\sinc^2(\xi_j),
 \qquad
 \sinc(x)=\frac{\sin\pi x}{\pi x}.
\end{equation}
For the cyclic families realizing the exponent \(1/n\), every fixed nonexceptional dual mode is quadratic in the cycle monomial, whereas the total energy is linear in its absolute value. Normalized spectral mass therefore escapes every finite frequency window.

The article is organized as follows. Sections~\ref{sec:overlap}--\ref{sec:global-locus} establish the global energy and branch geometry. Sections~\ref{sec:poset} and \ref{sec:scheme} develop the incidence poset and its local algebra. Section~\ref{sec:normal-form} proves the standard-lattice contact normal form. Sections~\ref{sec:valuations} and \ref{sec:fixed-window} classify local exponents and prove the sharp global estimate. Section~\ref{sec:planar} records the planar refinement, Section~\ref{sec:spectral} proves the primal--dual identity and spectral escape, and Section~\ref{sec:scope} proves the all-strata poset normal form.

\section{Overlap energy on lattice moduli}\label{sec:overlap}

Let
\[
 Q_n=[-1/2,1/2]^n,
 \qquad
 g_n(x)=|Q_n\cap(Q_n+x)|=\prod_{i=1}^n(1-|x_i|)_+.
\]
For a unimodular lattice \(\Lambda\subset\R^n\), define
\begin{equation}\label{eq:energy-def}
 \cE_n(\Lambda)=\sum_{v\in\Lambda\setminus\{0\}}g_n(v).
\end{equation}
The sum is finite because \(g_n\) is supported on \([-1,1]^n\).

Let
\[
 N_\Lambda(x)=\sum_{v\in\Lambda}\one_{Q_n}(x-v)
\]
be the occupation multiplicity on the torus \(\R^n/\Lambda\).

\begin{lemma}[variance identity]\label{lem:variance}
For every unimodular lattice,
\begin{equation}\label{eq:variance}
 \cE_n(\Lambda)
 =\int_{\R^n/\Lambda}(N_\Lambda-1)^2.
\end{equation}
\end{lemma}

\begin{proof}
Unfolding the square gives
\[
 \int_{\R^n/\Lambda}N_\Lambda^2
 =\sum_{v\in\Lambda}g_n(v)=1+\cE_n(\Lambda).
\]
Both the cube and the torus have volume one, so \(\int N_\Lambda=1\). Expanding the variance proves \eqref{eq:variance}.
\end{proof}

\begin{theorem}[proper overlap energy]\label{thm:proper}
The function \(\cE_n:X_n\to[0,\infty)\) is continuous and proper. More precisely,
\begin{equation}\label{eq:cusp-bound}
 \cE_n(\Lambda)\ge 2^{1-n}
 \left\lfloor\frac{1}{2\lambda_\infty(\Lambda)}\right\rfloor.
\end{equation}
\end{theorem}

\begin{proof}
Choose \(0\ne v\in\Lambda\) with \(\|v\|_\infty=s=\lambda_\infty(\Lambda)\), and put \(M=\lfloor(2s)^{-1}\rfloor\). For \(1\le m\le M\), every coordinate of \(mv\) has absolute value at most \(1/2\). Hence \(g_n(mv)\ge2^{-n}\). The \(2M\) vectors \(\pm mv\) are distinct, and \eqref{eq:cusp-bound} follows.

If \(\cE_n(\Lambda)\le R\), then
\[
 \lambda_\infty(\Lambda)>
 \frac{1}{2(2^{n-1}R+1)}.
\]
Mahler's compactness criterion \cite{Cassels,GruberLekkerkerker,Mahler} gives relative compactness of the sublevel. To prove continuity, choose a basis matrix \(B_0\). In a sufficiently small neighborhood of \(B_0\), only finitely many \(k\in\Z^n\) can satisfy \(Bk\in[-1,1]^n\). The local expression
\[
 B\longmapsto\sum_{0\ne k\in\Z^n}g_n(Bk)
\]
is therefore a finite sum of continuous functions. The sublevel is closed and hence compact.
\end{proof}

Let \(\cT_n\subset X_n\) denote the lattices for which \(Q_n+\Lambda\) tiles \(\R^n\) almost everywhere.

\begin{theorem}[zero set]\label{thm:zero-set}
For every unimodular lattice,
\[
 \cE_n(\Lambda)=0\quad\Longleftrightarrow\quad\Lambda\in\cT_n.
\]
In particular, \(\cT_n\) is compact.
\end{theorem}

\begin{proof}
The implication from tiling to zero energy follows from \cref{lem:variance}. Conversely, zero energy says that \(g_n(v)=0\) for every nonzero \(v\in\Lambda\). Thus the interiors of distinct cube translates are disjoint. The packing multiplicity on the quotient is at most one almost everywhere and has mean one, so it equals one almost everywhere. Compactness follows from \cref{thm:proper}.
\end{proof}

Fix any smooth Riemannian metric on \(X_n\). Different choices are locally bi-Lipschitz on compact sets.

\begin{theorem}[global qualitative error bound]\label{thm:global-holder}
There exist \(C_n,\theta_n,\delta_n>0\) such that
\begin{equation}\label{eq:global-holder}
 \dist_{X_n}(\Lambda,\cT_n)
 \le C_n\cE_n(\Lambda)^{\theta_n}
\end{equation}
whenever \(\cE_n(\Lambda)\le\delta_n\).
\end{theorem}

\begin{proof}
Choose a basis chart around \(\Lambda_0\in\cT_n\) on which the right action of \(\mathrm{SL}_n(\Z)\) is injective. Only finitely many deck labels contribute there, so the energy restricted to \(\det B=1\) is continuous and semialgebraic. Its zero set is the local lift of \(\cT_n\). The semialgebraic error-bound theorem \cite{BasuPollackRoy,BochnakCosteRoy} yields a local estimate with some positive exponent. The compactness of \(\cT_n\) supplies a finite chart cover. Finally, properness and \cref{thm:zero-set} imply that a sufficiently small energy sublevel lies in that cover.
\end{proof}

\section{The global Haj\'os locus}\label{sec:global-locus}

Let \(N_n(\R)\) be the group of upper unitriangular matrices and let \(N_n(\Z)=N_n(\R)\cap\mathrm{SL}_n(\Z)\). For each coordinate order \(\sigma\in S_n\), choose an orientation-preserving signed permutation matrix \(W_\sigma\) inducing that order, and set
\begin{equation}\label{eq:Ysigma}
 \cY_\sigma=
 \{W_\sigma^{-1}U\Z^n:U\in N_n(\R)\}\subset X_n.
\end{equation}

\begin{theorem}[nilmanifold branch decomposition]\label{thm:nil-branches}
For every \(n\ge2\),
\begin{equation}\label{eq:nil-union}
 \cT_n=\bigcup_{\sigma\in S_n}\cY_\sigma.
\end{equation}
Each \(\cY_\sigma\) is an embedded compact nilmanifold canonically diffeomorphic to \(N_n(\R)/N_n(\Z)\), of dimension \(\binom n2\).
\end{theorem}

\begin{proof}
The lattice form of the Minkowski--Haj\'os theorem says that a determinant-one basis matrix \(B\) generates a cube tiling if and only if there exist a coordinate permutation \(P\) and \(R\in\mathrm{SL}_n(\Z)\) for which \(PBR\) is unitriangular; see \cite{Hajos,Kolountzakis1998,LeePfanderWalnut}. Orientation signs may be moved between \(P\) and \(R\), giving \eqref{eq:nil-union}.

Two matrices \(U_1,U_2\in N_n(\R)\) generate the same lattice precisely when \(U_2^{-1}U_1\in N_n(\Z)\). The induced map from \(N_n(\R)/N_n(\Z)\) to \(X_n\) is therefore injective and immersive. The integer unitriangular group is a cocompact lattice in the simply connected nilpotent group \(N_n(\R)\) \cite{Raghunathan}; hence the source is compact. A proper injective immersion into the Hausdorff manifold \(X_n\) is an embedding.
\end{proof}

For \(\Lambda_0=B_0\Z^n\in\cT_n\), let \(\mathcal F(\Lambda_0)\) be the coordinate orders for which \(\Lambda_0\in\cY_\sigma\).

\begin{proposition}[local affine atlas]\label{prop:affine-atlas}
There is an injective basis chart \(\mathcal V\subset\mathrm{SL}_n(\R)\) about \(B_0\) and, for each \(\sigma\in\mathcal F(\Lambda_0)\), a fixed \(R_\sigma\in\mathrm{SL}_n(\Z)\) such that
\begin{equation}\label{eq:affine-atlas}
 \{B\in\mathcal V:B\Z^n\in\cT_n\}
 =\bigcup_{\sigma\in\mathcal F(\Lambda_0)}
 \{B\in\mathcal V:W_\sigma BR_\sigma\in N_n(\R)\}.
\end{equation}
Every member of the union is the intersection of \(\mathcal V\) with an affine subspace of matrix space.
\end{proposition}

\begin{proof}
Proper discontinuity of the integral basis action gives an injective chart. For each incident order, choose \(R_\sigma\) at \(B_0\). 
By discreteness, reduce the chart so that the same lift is the only one meeting it, modulo right multiplication by $N_n(\Z)$. Such right multiplication preserves the unitriangular condition.
The nonincident compact branches miss \(B_0\), so the chart may also be chosen disjoint from them. The equations fixing the diagonal to one and the entries below the chosen order to zero are affine linear.
\end{proof}

\begin{remark}
The standard lattice belongs to all \(n!\) branches. A generic point on any one branch belongs to no other branch. The following section gives a canonical description of every intermediate incidence.
\end{remark}

\section{The canonical incidence poset}\label{sec:poset}

For \(S\subset[n]\), put
\[
 E_S=\operatorname{span}\{e_i:i\in S\},
 \qquad
 \Lambda_S=\Lambda\cap E_S.
\]
The set \(S\) is called \emph{tight} if \(\Lambda_S\) has full rank in \(E_S\) and covolume one. The empty set is tight by convention. A rank-deficient section is assigned infinite covolume.

For a \(\Lambda\)-rational subspace \(V\), write
\[
 H_\Lambda(V)=\covol(\Lambda\cap V).
\]

\begin{lemma}[height submodularity]\label{lem:height-submodular}
If \(V\) and \(W\) are \(\Lambda\)-rational, then
\begin{equation}\label{eq:height-submodular}
 H_\Lambda(V+W)H_\Lambda(V\cap W)
 \le H_\Lambda(V)H_\Lambda(W).
\end{equation}
\end{lemma}

\begin{proof}
Let \(L_0=\Lambda\cap V\cap W\), \(L_V=\Lambda\cap V\), and \(L_W=\Lambda\cap W\). Choose a basis \(c\) of \(L_0\) and extend it to bases \((c,a)\) of \(L_V\) and \((c,b)\) of \(L_W\). The exterior-algebra volume inequality gives
\[
 \|c\wedge a\|\,\|c\wedge b\|
 \ge \|c\|\,\|c\wedge a\wedge b\|.
\]
The lattice generated by \(c,a,b\) has finite index at least one in \(\Lambda\cap(V+W)\). Consequently
\[
 \|c\wedge a\wedge b\|\ge H_\Lambda(V+W),
\]
which proves \eqref{eq:height-submodular}. Related determinant identities for rational sections and projections appear in \cite{McMullen}.
\end{proof}

\begin{lemma}[tight sets and flag prefixes]\label{lem:tight-prefix}
A coordinate set \(S\) is tight if and only if it occurs as a prefix of a unitriangular Haj\'os flag of \(\Lambda\).
Moreover, a maximal coordinate chain
\[
 \varnothing=S_0\subset S_1\subset\cdots\subset S_n=[n],
 \qquad |S_k|=k,
\]
is the prefix chain of a single Haj\'os flag if and only if every \(S_k\) is tight.
\end{lemma}

\begin{proof}
A flag prefix has a unitriangular basis and therefore covolume one. Conversely, suppose that \(S\) is tight. Since \(Q_n+\Lambda\) is a packing, \(Q_S+\Lambda_S\) is a packing in \(E_S\). Its tile and lattice covolumes are both one, so it is a tiling.

The section lattice is primitive in \(\Lambda\). Projection to \(E_{S^c}\) gives a lattice
\[
 \overline\Lambda=\operatorname{pr}_{S^c}(\Lambda)
 \simeq\Lambda/\Lambda_S
\]
of covolume one. 
It tiles by $Q_{S^c}$. To see this, periodize the product cube and group the sum over $\Lambda$ into cosets of $\Lambda_S$. In each coset sum, every point of $E_S$ is covered exactly once except possibly on the boundaries of the translated copies of $Q_S$. The union of these boundaries has Lebesgue measure zero in $E_S$, so the inner sum equals one almost everywhere. The remaining sum is the $Q_{S^c}$-periodization.
Apply the Minkowski--Haj\'os theorem to the section and quotient, lift a unitriangular quotient basis, and append it to a unitriangular section basis. 
In the coordinate order listing \(S\) first, the resulting basis is block upper unitriangular.

For the final assertion, suppose every member of the maximal chain is tight and write \(S_k=S_{k-1}\cup\{i_k\}\). Put \(L_k=\Lambda\cap E_{S_k}\). The primitive quotient \(L_k/L_{k-1}\) is infinite cyclic. Choose \(v_k\in L_k\) whose class generates it. If a basis of \(L_{k-1}\) is followed by \(v_k\), the covolume ratio is \(|(v_k)_{i_k}|\); since both covolumes equal one, this coordinate is \(\pm1\). Reverse \(v_k\) when necessary. The vectors \(v_1,\ldots,v_n\) form a basis of \(\Lambda\), and their matrix is upper unitriangular in the order \(i_1,\ldots,i_n\). The converse is immediate from the first assertion.
\end{proof}

\begin{theorem}[incidence-poset theorem]\label{thm:incidence-poset}
The tight coordinate subsets of a lattice cube tiler form a sublattice of the Boolean lattice. 
There is a unique partial order \(P_\Lambda\) on \([n]\) such that
\begin{equation}\label{eq:tight-ideals}
 S\text{ is tight}
 \quad\Longleftrightarrow\quad
 S\text{ is an order ideal of }P_\Lambda.
\end{equation}
The admissible Haj\'os flags of \(\Lambda\) are exactly the linear extensions of \(P_\Lambda\).
\end{theorem}

\begin{proof}
Every full-rank coordinate section of a cube-packing lattice has covolume at least one. If \(S\) and \(T\) are tight, \cref{lem:height-submodular} gives
\[
 H_\Lambda(E_{S\cup T})H_\Lambda(E_{S\cap T})\le1.
\]
Both factors are at least one, hence both equal one. Thus tight sets are closed under union and intersection.

Define
\begin{equation}\label{eq:poset-def}
 i\preceq_\Lambda j
 \quad\Longleftrightarrow\quad
 i\in S\text{ for every tight }S\text{ containing }j.
\end{equation}
This is a preorder. 
A full Haj\'os flag supplies a chain of tight sets that separates every pair of coordinates, so the preorder is antisymmetric. Every tight set is a down-set. Conversely, the smallest tight set containing \(j\) is the intersection of all tight sets containing it. A down-set is the union of these smallest tight sets over its elements and is therefore tight. This proves \eqref{eq:tight-ideals}. Maximal chains of order ideals are the linear extensions of the poset \cite{StanleyEC}; \cref{lem:tight-prefix} identifies them with admissible flags.
\end{proof}

Fix a tiler \(\Lambda=U\Z^n\) and use the right-invariant coordinate
\begin{equation}\label{eq:right-coordinate}
 B=(I+A)U.
\end{equation}
For an incident order \(\sigma\), choose \(R_\sigma\in\mathrm{SL}_n(\Z)\) such that \(W_\sigma UR_\sigma\) is upper unitriangular. Then
\[
 W_\sigma(B-U)R_\sigma
 =(W_\sigma A W_\sigma^{-1})(W_\sigma UR_\sigma).
\]
Right multiplication by an upper unitriangular matrix preserves and reflects strict upper triangularity. 
Consequently the local branches become simultaneous coordinate subspaces:
\begin{equation}\label{eq:poset-arrangement}
 \cH(P_\Lambda)
 =\bigcup_{\sigma\in\LinExt(P_\Lambda)}
 \{A:a_{ii}=0,\ a_{ij}=0\text{ whenever }j<_\sigma i\}.
\end{equation}

\begin{theorem}[universality]\label{thm:universality}
Every finite poset \(P\) on \([n]\) occurs as \(P_\Lambda\) for a lattice cube tiler. More explicitly, choose a natural labeling of \(P\), and for every \(i<_Pj\) choose \(\alpha_{ij}\in(0,1)\) so that
\begin{equation}\label{eq:row-independence}
 \{1\}\cup\{\alpha_{ij}:i<_Pj\}
 \quad\text{is linearly independent over }\mathbb Q
\end{equation}
for each fixed row \(i\). Define
\begin{equation}\label{eq:UP}
 (U_P)_{ii}=1,
 \qquad
 (U_P)_{ij}=\begin{cases}
 \alpha_{ij},&i<_Pj,\\
 0,&\text{otherwise}.
 \end{cases}
\end{equation}
Then \(U_P\Z^n\) tiles by \(Q_n\) and \(P_{U_P\Z^n}=P\).
\end{theorem}

\begin{proof}
The matrix \(U_P\) is upper unitriangular, so its lattice tiles. If \(S\) is an order ideal, every column indexed by \(S\) is supported in \(E_S\), and the corresponding principal submatrix is unitriangular. Thus \(S\) is tight.

Conversely, let \(U_Pm\in E_S\), with \(m\in\Z^n\). For \(i\notin S\), the \(i\)-th coordinate equation is
\[
 m_i+\sum_{j:\,i<_Pj}\alpha_{ij}m_j=0.
\]
The rowwise independence \eqref{eq:row-independence} forces \(m_i=0\) and \(m_j=0\) for every successor \(j\) of \(i\). Therefore \(m\) is supported on the largest order ideal contained in \(S\),
\[
 I(S)=\{j\in S:\downarrow_Pj\subset S\}.
\]
If \(S\) is not an ideal, then \(|I(S)|<|S|\), and the coordinate section is rank deficient. The characterization \eqref{eq:tight-ideals} now gives \(P_{U_P\Z^n}=P\). The parameters may be chosen algebraic, for instance by using sufficiently small rational multiples of square roots of distinct primes in each row.
\end{proof}

\begin{corollary}\label{cor:branch-complexity}
Every finite linear-extension pattern occurs as the branch incidence of a lattice cube tiler. In particular, the number of local irreducible branches can realize the \(\#P\)-complete function \(P\mapsto|\LinExt(P)|\).
\end{corollary}

\section{Rooted circuits and local scheme invariants}\label{sec:scheme}

Regard every relation \(i<_Pj\) as a fixed directed edge \(i\to j\). The off-diagonal coordinates \(a_{ij}\) supply variable directed edges. A \emph{rooted \(P\)-circuit} is the variable-edge set of an inclusion-minimal directed cycle in the union of \(P\) with the variable support. This terminology is compatible with the standard circuit language for oriented combinatorial structures \cite{BjornerEtAl}.

Let \(R=\C[x_{ij}:i\ne j]\), and define
\begin{equation}\label{eq:rooted-ideal}
 I_P=\left\langle
 \prod_{(i,j)\in C}x_{ij}:C\text{ is a rooted }P\text{-circuit}
 \right\rangle.
\end{equation}

\begin{theorem}[rooted-circuit ideal]\label{thm:rooted-ideal}
The ideal \(I_P\) is the squarefree Stanley--Reisner ideal of the off-diagonal branch arrangement in \eqref{eq:poset-arrangement}. More precisely,
\begin{equation}\label{eq:prime-decomp}
 I_P=\bigcap_{\sigma\in\LinExt(P)}\mathfrak p_\sigma,
 \qquad
 \mathfrak p_\sigma=\langle x_{ij}:j<_\sigma i\rangle.
\end{equation}
It is radical and equidimensional of height \(\binom n2\).
\end{theorem}

\begin{proof}
An edge support lies in a branch precisely when all its edges point forward in some linear extension of \(P\). Equivalently, its union with \(P\) is acyclic. The minimal forbidden supports are exactly the rooted circuits. The maximal allowed supports are the forward tournaments
\[
 F_\sigma=\{(i,j):i<_\sigma j\},
 \qquad \sigma\in\LinExt(P).
\]
Each has \(\binom n2\) elements. The Stanley--Reisner facet--prime correspondence \cite{MillerSturmfels} proves the result.
\end{proof}

For a real matrix \(A\), put
\begin{equation}\label{eq:rhoP}
 \rho_P(A)=\max_C\min_{e\in C}|a_e|,
\end{equation}
with value zero if there is no rooted circuit in the support.

\begin{theorem}[bottleneck duality]\label{thm:poset-bottleneck}
For every \(A\),
\begin{align}
 \rho_P(A)
 &=\inf\{t\ge0:P\cup\{e:|a_e|>t\}\text{ is acyclic}\},\label{eq:rho-threshold}\\
 &=\min_{\sigma\in\LinExt(P)}\max_{j<_\sigma i}|a_{ij}|.\label{eq:rho-order}
\end{align}
Consequently,
\begin{equation}\label{eq:exact-branch-distance}
 \dist_\infty(A,\cH(P))
 =\max\left\{\max_i|a_{ii}|,\rho_P(A)\right\}.
\end{equation}
\end{theorem}

\begin{proof}
The threshold graph is acyclic exactly when every rooted circuit contains an edge of magnitude at most \(t\), proving \eqref{eq:rho-threshold}. An acyclic extension of \(P\) has a topological order, and every such order is a linear extension of \(P\), which proves \eqref{eq:rho-order}. For a fixed order, the closest branch point keeps every forward entry, deletes every backward entry, and sets the diagonal to zero. Minimizing over the orders gives \eqref{eq:exact-branch-distance}.
\end{proof}

Let
\[
 c(P)=\max\{|C|:C\text{ is a rooted }P\text{-circuit}\}.
\]
Every rooted circuit is carried by a simple directed cycle and has distinct tails; hence \(c(P)\le n\).

\begin{corollary}[algebraic repair exponent]\label{cor:algebraic-repair}
In a unit neighborhood of the origin,
\begin{equation}\label{eq:algebraic-repair}
 \dist_\infty(A,\cH(P))
 \le C_P
 \left(
   \sum_i|a_{ii}|+
   \sum_C\prod_{e\in C}|a_e|
 \right)^{1/c(P)}.
\end{equation}
The exponent \(1/c(P)\) is optimal for the rooted-circuit gauge.
\end{corollary}

\begin{proof}
Let \(\rho=\rho_P(A)\). A circuit attaining the maximum in \eqref{eq:rhoP} has \(\ell\le c(P)\) edges, each of magnitude at least \(\rho\). Since \(\rho\le1\), its monomial is at least \(\rho^\ell\ge\rho^{c(P)}\). Combine this with \eqref{eq:exact-branch-distance}. For optimality, put the same parameter on the edges of a rooted circuit of maximum size and set the other normal coordinates to zero.
\end{proof}

\begin{theorem}[multiplicity and log-canonical threshold]\label{thm:scheme-invariants}
Let \(e(P)=|\LinExt(P)|\). Then
\begin{equation}\label{eq:multiplicity}
 \mult_0(I_P)=e(P),
 \qquad
 \lct_0(I_P)=\binom n2.
\end{equation}
Inside the smooth determinant-one hypersurface, let \(J_P\) be the ideal obtained by adjoining the local diagonal equations. Then
\begin{equation}\label{eq:lct-slice}
 \lct_0(J_P)=\binom n2+n-1
 =\operatorname{codim}V(J_P).
\end{equation}
\end{theorem}

\begin{proof}
The prime decomposition \eqref{eq:prime-decomp} is a reduced equidimensional union of \(e(P)\) coordinate planes, each of degree one, and distinct components meet in smaller dimension. Hence the degree, equivalently the Hilbert--Samuel multiplicity at the origin, is \(e(P)\).

Howald's formula \cite{Howald} identifies the log-canonical threshold with the fractional circuit-cover program
\begin{equation}\label{eq:lct-lp}
 \min\left\{
 \sum_{i\ne j}u_{ij}:u_{ij}\ge0,\
 \sum_{(i,j)\in C}u_{ij}\ge1
 \text{ for every rooted circuit }C
 \right\}.
\end{equation}
For every comparable pair \(i<_Pj\), the reverse edge \(j\to i\) is a singleton circuit, so \(u_{ji}\ge1\). For every incomparable pair, the directed two-cycle gives \(u_{ij}+u_{ji}\ge1\). Summing these disjoint constraints over the \(\binom n2\) unordered pairs gives the lower bound \(\binom n2\).

For the reverse inequality, assign weights as follows: (i) weight one to each edge opposite a comparison; (ii) weight zero to each edge parallel to a comparison; and (iii) weight $1/2$ to each orientation of an incomparable pair. If \(i<_Pj\), a minimal rooted circuit cannot use the variable edge \(i\to j\): replacing that occurrence by the fixed \(P\)-edge \(i\to j\), or equivalently by a fixed directed \(P\)-path, produces the same directed cycle with a strictly smaller variable support. This contradicts inclusion-minimality.

It remains to check the weights on a rooted circuit. If the circuit contains an edge opposite a comparison, that edge alone has weight one. Otherwise every variable edge joins an incomparable pair. There cannot be exactly one such edge \(i\to j\), since the rest of the directed cycle would then give a fixed \(P\)-path from \(j\) to \(i\), making the pair comparable. Hence there are at least two incomparable-pair edges, each of weight \(1/2\). Every circuit constraint in \eqref{eq:lct-lp} is therefore satisfied. For each unordered pair the two oriented weights sum to one, so the total weight is \(\binom n2\).

On the determinant-one slice, eliminate one diagonal coordinate. Modulo the other \(n-1\) diagonal coordinates and \(I_P\), every nonidentity determinant monomial vanishes because it contains a directed cycle. The eliminated diagonal therefore lies in the same local ideal. Thus \(J_P\) is analytically the sum of \(I_P\) and \(n-1\) independent coordinate variables. Adding their singleton constraints to \eqref{eq:lct-lp} proves \eqref{eq:lct-slice}.
\end{proof}

Let \(\mathsf{Acyc}_n\) be the simplicial complex of acyclic directed graphs on \([n]\). After the cone variables parallel to relations of \(P\) are removed, the Stanley--Reisner complex of \(I_P\) is
\begin{equation}\label{eq:link}
 \Delta_P^{\mathrm{ess}}=\link_{\mathsf{Acyc}_n}(P).
\end{equation}
At the standard lattice \(P\) is the antichain, and the theorem of Bj\"orner and Welker gives
\begin{equation}\label{eq:acyc-sphere}
 \mathsf{Acyc}_n\simeq S^{n-2}.
\end{equation}
The broader topology of directed graph complexes is treated in \cite{BjornerWelker}.

\section{The standard-lattice contact normal form}\label{sec:normal-form}

For \(B=(b_{ij})\in\mathrm{SL}_n(\R)\), put
\begin{equation}\label{eq:DeltaB}
 \Delta_n(B)=\sum_{k\in\Z^n\setminus\{0\}}g_n(Bk).
\end{equation}
Let
\begin{equation}\label{eq:Hn}
 \cH_n=\bigcup_{\sigma\in S_n}P_\sigma N_n(\R)P_\sigma^{-1}.
\end{equation}
Thus \(\cH_n\) is the union of the unitriangular branches through the identity.

\begin{lemma}[local zero locus]\label{lem:local-zero}
There is a neighborhood \(\mathcal U\) of \(I\) in \(\mathrm{SL}_n(\R)\) such that
\[
 \Delta_n(B)=0\quad\Longleftrightarrow\quad B\in\cH_n
\]
for every \(B\in\mathcal U\).
\end{lemma}

\begin{proof}
Every matrix in \(\cH_n\) generates a cube-tiling lattice. Conversely, zero overlap gives a unimodular cube packing and therefore a tiling. Haj\'os's theorem gives \(e_i\in B\Z^n\) for some \(i\). Reduce the neighborhood so that \(\|B^{-1}-I\|_{\infty\to\infty}<1/2\). If \(e_i=Bk\), then the integer vector \(k=B^{-1}e_i\) is componentwise within \(1/2\) of \(e_i\), so \(k=e_i\). The \(i\)-th column of \(B\) is therefore \(e_i\).

Delete row and column \(i\). The resulting matrix \(C\) has determinant one. Its cube translates form a packing: otherwise a nonzero projected lattice vector in \((-1,1)^{n-1}\), after addition of a suitable integer multiple of \(e_i\), would give a nonzero vector of \(B\Z^n\) in \((-1,1)^n\). Hence \(C\) gives a lower-dimensional cube tiling. Induction puts \(C\), and then \(B\), in a unitriangular branch.
\end{proof}

For a directed simple cycle \(C\), counted modulo cyclic rotation, define
\[
 m_C(B)=\prod_{(i,j)\in C}b_{ij},
\]
and set
\begin{equation}\label{eq:augmented-gauge}
 \cA(B)=
 \sum_{i=1}^n|b_{ii}-1|+\sum_C|m_C(B)|.
\end{equation}

The proof of the normal form uses convergent real Puiseux arcs. After a finite reparametrization, this is equivalent to the semialgebraic arc language used below.

\begin{lemma}[minimum-cycle row isolation]\label{lem:row-isolation}
Give every off-diagonal edge \(i\to j\) a positive weight \(w_{ij}\), allowing \(+\infty\). Let \(C=(i_1,\ldots,i_\ell)\) be a directed simple cycle of minimum finite total weight \(\mu\). For a vertex \(i_r\in C\), every chord \(i_r\to i_s\), \(i_s\ne i_{r+1}\), satisfies
\[
 w_{i_ri_s}>w_{i_ri_{r+1}}.
\]
Two distinct minimum-weight cycles cannot have the same vertex set.
\end{lemma}

\begin{proof}
Suppose a chord has weight no larger than the cycle edge leaving its tail. Replace that edge by the chord and retain the successor edge at every other vertex of \(C\). The resulting functional digraph contains a proper directed cycle and at least one tree edge feeding it. Its total selected weight is at most \(\mu\), while the proper cycle has strictly smaller weight. This contradicts minimality. If two distinct minimum cycles had the same vertices, one of the changed rows would use a chord of no greater weight, giving the same contradiction.
\end{proof}

Let \(B(t)\in\mathrm{SL}_n(\R)\) be a nonconstant convergent Puiseux arc with \(B(0)=I\). For \(i\ne j\), write
\[
 b_{ij}(t)=a_{ij}t^{w_{ij}}+o(t^{w_{ij}}),
 \qquad a_{ij}\ne0,
\]
and set \(w_{ij}=+\infty\) when the entry vanishes identically. Define
\begin{equation}\label{eq:lambda-mu}
 \lambda=\min_i\ord_t(b_{ii}-1),
 \qquad
 \mu=\min_C\sum_{(i,j)\in C}w_{ij},
 \qquad
 \nu=\min\{\lambda,\mu\}.
\end{equation}
At order \(\nu\), write
\[
 b_{ii}(t)=1+d_it^\nu+o(t^\nu),
\]
with \(d_i=0\) if the diagonal order is larger. For each cycle of weight \(\nu\), put
\[
 a_C=\prod_{(i,j)\in C}a_{ij},
 \qquad
 \chi_C=(-1)^{|C|-1}a_C,
\]
and define
\begin{equation}\label{eq:PNR}
 P=\sum_C[\chi_C]_+,\qquad
 N=\sum_C[-\chi_C]_+,\qquad
 R=\sum_i[-d_i]_+.
\end{equation}
The cycle sums vanish when \(\mu>\nu\).

\begin{lemma}[first-contact classification]\label{lem:first-contact}
For every \(0\ne k\in\{-1,0,1\}^n\), either the summand \(g_n(B(t)k)\) vanishes identically near zero or
\[
 \ord_t g_n(B(t)k)\ge\nu.
\]
Equality is possible only in one of the following cases:
\begin{enumerate}[label=\textup{(\roman*)}]
\item \(k=\pm e_i\) and \(d_i<0\);
\item the support of \(k\) is the vertex set of a minimum-weight directed cycle \(C\) with \(\chi_C<0\), and the signs of \(k\) are the compatible cycle signs, unique up to overall reversal.
\end{enumerate}
\end{lemma}

\begin{proof}
Let \(S=\{i:k_i\ne0\}\). In each row \(i\in S\), a positive overlap factor requires the perturbation of the base value \(k_i\) to point inward. Its order is at least the minimum of the diagonal-deviation order in that row and the weights \(w_{ij}\) with \(j\in S\setminus\{i\}\). Select one term attaining this lower row order. A diagonal choice is represented by a self-loop; an off-diagonal choice by an edge \(i\to j\) on \(S\).

If a self-loop is selected, the sum of the row orders is at least \(\lambda\). If every choice is off diagonal, the selected functional digraph contains a directed cycle, so the sum is at least \(\mu\). Cancellation inside a row can only increase its actual order. This proves the lower bound.

If equality comes from a self-loop, every additional row would add positive order, so \(S\) is a singleton. If equality comes from off-diagonal choices, there is no tree vertex outside the selected cycle and that cycle has weight \(\mu=\nu\). Equality also excludes cancellation in every selected row. By \cref{lem:row-isolation}, the cycle edge is the unique leading outgoing edge in its row. The inward sign equations are
\[
 k_{i_{r+1}}=-\operatorname{sgn}(a_{i_ri_{r+1}})k_{i_r}.
\]
They close exactly when \(\operatorname{sgn}(a_C)=(-1)^{|C|}\), which is equivalent to \(\chi_C<0\).
\end{proof}

\begin{theorem}[exact tropical contact]\label{thm:tropical-contact}
Unless \(B(t)\) lies identically in \(\cH_n\),
\begin{equation}\label{eq:arc-orders}
 \ord_t\Delta_n(B(t))
 =\nu
 =\ord_t\cA(B(t)).
\end{equation}
The two functions have the same leading coefficient:
\begin{equation}\label{eq:leading-coefficient}
 [t^\nu]\Delta_n(B(t))
 =2(N+R)
 =[t^\nu]\cA(B(t)).
\end{equation}
\end{theorem}

\begin{proof}
After reducing the matrix neighborhood, a label with \(\max_i|k_i|\ge2\) has some coordinate of \(Bk\) outside \([-1,1]\), so only \(\{-1,0,1\}^n\) contributes. By \cref{lem:first-contact}, the two compatible labels of each negative-sign minimum cycle contribute \(2|a_C|t^\nu\), and the two singleton labels associated with each \(d_i<0\) contribute \(2|d_i|t^\nu\). Minimum cycles with the same support are excluded by \cref{lem:row-isolation}. Thus the coefficient of the energy is \(2(N+R)\).

At order \(\nu\), the determinant expansion contains only the diagonal term and the single minimum-cycle permutation terms. Products of two nontrivial cycles have order greater than \(\nu\). Hence
\begin{equation}\label{eq:det-balance}
 \sum_i d_i+P-N=0.
\end{equation}
If \(S=\sum_i[d_i]_+\), then \(S-R=N-P\), and the coefficient of the gauge is
\[
 S+R+P+N=2(N+R).
\]
This is positive unless all order-\(\nu\) diagonal and cycle terms vanish, contrary to the definition of \(\nu\). 
The order and coefficient identities follow.
\end{proof}

\begin{theorem}[constant-one normal form]\label{thm:normal-form}
For every \(n\ge2\),
\begin{equation}\label{eq:normal-form}
 \lim_{\substack{B\to I,\ B\in\mathrm{SL}_n(\R)\\B\notin\cH_n}}
 \frac{\Delta_n(B)}{\cA(B)}=1.
\end{equation}
Equivalently,
\[
 \Delta_n(B)=(1+o(1))
 \left(
 \sum_i|b_{ii}-1|+\sum_C|m_C(B)|
 \right)
\]
uniformly near \(I\) on the determinant-one hypersurface.
\end{theorem}

\begin{proof}
Both functions are continuous and semialgebraic near \(I\), and \cref{lem:local-zero} identifies their common zero set. If \eqref{eq:normal-form} failed, there would be \(\varepsilon>0\) such that the semialgebraic set
\[
 \left\{B\notin\cH_n:
 \left|\frac{\Delta_n(B)}{\cA(B)}-1\right|\ge\varepsilon\right\}
\]
accumulated at \(I\). Semialgebraic curve selection \cite{BasuPollackRoy,BochnakCosteRoy} would then give a convergent Puiseux arc in this failure set, contradicting \cref{thm:tropical-contact}.
\end{proof}

Define
\begin{equation}\label{eq:rho-standard}
 \rho(B)=\max_C\min_{(i,j)\in C}|b_{ij}|,
\end{equation}
with value zero when the off-diagonal support is acyclic.

\begin{proposition}[exact bottleneck distance]\label{prop:standard-bottleneck}
For every real matrix \(B\),
\begin{equation}\label{eq:standard-bottleneck}
 \dist_\infty(B,\cH_n)
 =\max\left\{\max_i|b_{ii}-1|,\rho(B)\right\}.
\end{equation}
\end{proposition}

\begin{proof}
The graph formed by the edges of magnitude greater than \(t\) is acyclic exactly when every directed cycle has an edge of magnitude at most \(t\). Thus \(\rho(B)\) is the least threshold producing an acyclic graph. A topological order of that graph selects a Haj\'os branch. For a fixed order, the nearest branch point keeps the forward entries, deletes the backward entries, and sets the diagonal to one. Minimization over all orders proves the formula.
\end{proof}

\begin{theorem}[sharp standard-lattice rigidity]\label{thm:sharp-standard}
There are \(K_n,\varepsilon_n>0\) such that
\begin{equation}\label{eq:sharp-standard}
 \dist_\infty(B,\cH_n)
 \le K_n\Delta_n(B)^{1/n}
\end{equation}
whenever \(B\in\mathrm{SL}_n(\R)\) and \(\|B-I\|_\infty<\varepsilon_n\). The exponent \(1/n\) is optimal.
\end{theorem}

\begin{proof}
Let \(\rho=\rho(B)\). A cycle attaining \eqref{eq:rho-standard} has length \(\ell\le n\), and all its edges have magnitude at least \(\rho\). In a unit neighborhood,
\[
 \sum_C|m_C(B)|\ge\rho^\ell\ge\rho^n.
\]
The normal form also gives \(\max_i|b_{ii}-1|\le2\Delta_n(B)\) after reducing the neighborhood. Apply \cref{prop:standard-bottleneck}.

For optimality, put the same parameter \(t\) on the edges of a directed Hamiltonian cycle and use the determinant-one correction of \cref{eq:cyclic-family} below. The distance to every acyclic branch is at least \(|t|\), while the energy is asymptotic to \(2|t|^n\).
\end{proof}

For a fixed directed graph \(G\), let \(\mathcal S_G\) be the determinant-one matrices near \(I\) whose off-diagonal support is contained in \(G\). Let \(c(G)\) be the maximum length of a directed simple cycle, with \(c(G)=1\) when \(G\) is acyclic.

\begin{theorem}[directed-circumference law]\label{thm:circumference}
The optimal uniform local repair exponent on \(\mathcal S_G\) is
\[
 \alpha_G=\frac1{c(G)}.
\]
\end{theorem}

\begin{proof}
The proof of \cref{thm:sharp-standard} gives \(\rho^{c(G)}\le2\Delta_n\), because every cycle realizing the bottleneck has length at most \(c(G)\). This proves the upper estimate. For sharpness, choose a directed cycle in \(G\) of length \(c(G)\), put the same parameter on its edges, and leave the complementary coordinates fixed. Its energy has order \(c(G)\) and its branch distance has order one. If \(G\) is acyclic, only the diagonal defect is normal to \(\cH_n\), and \cref{thm:normal-form} gives a linear estimate.
\end{proof}

\section{Signed contact valuations at arbitrary tilers}\label{sec:valuations}

The standard-lattice normal form uses the fact that all contact forms are coordinate perturbations. At a general tiler, several words may have the same order inside one affine boundary form and may cancel before the positive part is taken. The optimal local exponent is therefore attached to the complete signed contact products rather than to the reduced rooted-circuit ideal alone.

Fix a tiler \(U\Z^n\), write a nearby basis as
\[
 B=(I+A)U,\qquad \det(I+A)=1,
\]
and denote by \(\cH_U\) the finite local union of incident affine Haj\'os branches in the \(A\)-coordinate. After reducing the matrix neighborhood, one finite set of labels contains every possible contributor to \(\Delta_n(B)\). Subdivide the determinant-one germ by the signs of all coordinate slacks
\[
 1\pm((I+A)Uk)_i
\]
and by the Voronoi regions of the incident branches. On the relative interior of a resulting cell \(\Gamma\), one branch \(L_\Gamma\) realizes the distance to \(\cH_U\), up to equivalence of fixed norms, and
\begin{equation}\label{eq:contact-cell-sum}
 \Delta_n((I+A)U)=\sum_{a=1}^{s_\Gamma}p_{\Gamma,a}(A),
 \qquad p_{\Gamma,a}\ge0.
\end{equation}
Each \(p_{\Gamma,a}\) is a product of at most \(n\) complete affine slacks. Let \(J_\Gamma\) be the analytic ideal generated by these products and let \(N_\Gamma\) be the ideal generated by affine normal equations for \(L_\Gamma\). Then, on the cell,
\begin{equation}\label{eq:contact-ideal-comparison}
 \Delta_n((I+A)U)\asymp\|J_\Gamma(A)\|,
 \qquad
 \dist(A,\cH_U)\asymp\|N_\Gamma(A)\|.
\end{equation}

Take a common principalization of \(J_\Gamma N_\Gamma\) and the boundary ideal of \(\Gamma\). A prime divisor \(E\) on a toroidal refinement is called \emph{accessible} if its real locus meets the closure of the lifted relative interior of \(\Gamma\). 
After the sector-adapted toroidal refinement used below, an accessible divisor corresponds precisely to a divisorial ray contained in the closure of a lifted real sector cone associated with $\Gamma$; thus a monomial real arc in that sector can be chosen transverse to the divisor.
Write \(v_E\) for its valuation and set
\begin{equation}\label{eq:kappa-U}
 \kappa(\Gamma)=
 \max_{E\ {
m accessible}}
 \frac{v_E(J_\Gamma)}{v_E(N_\Gamma)},
 \qquad
 \kappa(U)=\max_\Gamma\kappa(\Gamma),
\end{equation}
where divisors with \(v_E(N_\Gamma)=0\) are omitted.

For a monomial arc transverse to \(E\), the valuation \(v_E(J_\Gamma)\) is the vanishing order of the contact ideal and hence, by \eqref{eq:contact-ideal-comparison}, of overlap on the cell. Similarly, \(v_E(N_\Gamma)\) is the vanishing order of the normal ideal and therefore of the distance to the local branch arrangement. Their quotient is precisely the contact order detected by that divisorial degeneration. The accessibility condition is essential here: a divisor on the complex resolution can affect the real local exponent only when its real locus is approached from the lifted real sign cell under consideration. Divisors visible only from other real sectors do not correspond to admissible degenerations in \(\Gamma\).

\begin{theorem}[signed contact-valuation theorem]\label{thm:signed-valuations}
The number \(\kappa(U)\) is a finite positive rational number. There are \(C_U,\varepsilon_U>0\) such that
\begin{equation}\label{eq:valuation-holder}
 \dist(A,\cH_U)
 \le C_U\Delta_n((I+A)U)^{1/\kappa(U)}
\end{equation}
on the determinant-one slice whenever \(\|A\|<\varepsilon_U\). The exponent is optimal. More precisely,
\begin{equation}\label{eq:arc-kappa}
 \kappa(U)=
 \max_\gamma
 \frac{\ord_t\Delta_n((I+\gamma(t))U)}
 {\ord_t\dist(\gamma(t),\cH_U)},
\end{equation}
where \(\gamma\) ranges over non-tiling convergent real Puiseux arcs based at the origin.
\end{theorem}

\begin{proof}
Principalization in characteristic zero \cite{Hironaka} and a toroidal refinement give finitely many real charts with coordinates \(x=(x_1,\ldots,x_r)\), analytic units \(u,v\), and exponent vectors \(a,b\in\Z_{\ge0}^r\) such that
\[
 \pi^*J_\Gamma=(u x^a),
 \qquad
 \pi^*N_\Gamma=(v x^b).
\]
Resolving the cell boundary at the same time makes every lifted sign sector a finite union of rational polyhedral orthants. For a monomial arc \(x_i=c_it^{w_i}\) in such a sector,
\[
 \ord_t\|J_\Gamma\|=a\cdot w,
 \qquad
 \ord_t\|N_\Gamma\|=b\cdot w.
\]
The two ideals have the same zero set on every accessible sector by \cref{eq:contact-ideal-comparison,thm:zero-set,prop:affine-atlas}. 
This remains true on the sector closure. A zero-energy limit belongs to the tiling arrangement. Moreover, the closure of the Voronoi region for $L_\Gamma$ can meet that zero set only on $L_\Gamma$. The semialgebraic {\L}ojasiewicz inequality therefore excludes an accessible ray with \(b\cdot w=0<a\cdot w\). On each rational polyhedral cone, the ratio \((a\cdot w)/(b\cdot w)\) attains its maximum on an extremal ray. After refinement that ray is divisorial, so the maximum is finite and rational.

The inequalities
\[
 v_E(J_\Gamma)\le\kappa(\Gamma)v_E(N_\Gamma)
\]
give \(\|N_\Gamma\|^{\kappa(\Gamma)}\le C_\Gamma\|J_\Gamma\|\) on every resolved chart. Together with \cref{eq:contact-ideal-comparison} and the finiteness of the cell decomposition, this proves \cref{eq:valuation-holder}. A monomial arc transverse to a maximizing accessible divisor realizes the quotient in \cref{eq:kappa-U}; its projection is a real Puiseux arc in the prescribed cell. No larger H\"older exponent can hold along that arc, proving both optimality and \cref{eq:arc-kappa}.
\end{proof}

\section{Fixed-window multiplicity and the sharp global bound}\label{sec:fixed-window}

The valuation formula identifies the exact local obstruction but does not by itself bound it uniformly. The required bound follows from a positive-product principle in which every affine factor is kept intact.

\begin{lemma}[positive-product error bound]\label{lem:positive-product}
Let \(P\subset\R^m\) be a compact polyhedron. Let \(\ell_h:P\to[0,\infty)\), \(h\in V\), be finitely many affine functions, and let \(\mathcal G\) be a finite family of nonempty subsets of \(V\), each of cardinality at most \(d\). For positive constants \(c_E\), define
\[
 F(x)=\sum_{E\in\mathcal G}c_E\prod_{h\in E}\ell_h(x),
 \qquad Z=\{x\in P:F(x)=0\}.
\]
If \(Z\ne\varnothing\), then
\begin{equation}\label{eq:positive-product-bound}
 \dist(x,Z)^d\le C_PF(x)
\end{equation}
for every \(x\in P\).
\end{lemma}

\begin{proof}
The assertion is immediate if \(\mathcal G\) is empty. Assume otherwise. Fix \(z\in Z\) and put \(V_z=\{h:\ell_h(z)=0\}\). Every edge of \(\mathcal G\) meets \(V_z\). On a neighborhood \(W_z\) of \(z\), all factors outside \(V_z\) are bounded below by a constant \(a_z>0\). Moreover,
\begin{equation}\label{eq:zero-cover-union}
 Z\cap W_z
 =\bigcup_{C}
 \bigl(P\cap W_z\cap\{\ell_h=0:h\in C\}\bigr),
\end{equation}
where \(C\subset V_z\) ranges over the vertex covers of \(\mathcal G\). Every linear system in this union is feasible because it contains \(z\). Hoffman's error bound, applied together with the inequalities defining \(P\), gives
\begin{equation}\label{eq:hoffman-local}
 \dist(x,Z)\le H_z\delta_z(x),
 \qquad
 \delta_z(x)=\min_{C\ {
m cover}}\max_{h\in C}\ell_h(x)
\end{equation}
on a smaller neighborhood.

For \(x\in Z\) the desired inequality is trivial. If \(x\notin Z\), then \(\delta_z(x)>0\), and there is an edge \(E_x\in\mathcal G\) such that
\begin{equation}\label{eq:blocker-edge}
 \ell_h(x)\ge\delta_z(x)
 \quad\text{for every }h\in E_x\cap V_z.
\end{equation}
Indeed, otherwise the vertices \(h\in V_z\) satisfying \(\ell_h(x)<\delta_z(x)\) would form a cover with smaller maximum. Since \(|E_x|\le d\), the unit lower bounds outside \(V_z\), \cref{eq:blocker-edge}, and \(\delta_z\le1\) imply
\[
 F(x)\ge c_z\delta_z(x)^d.
\]
Together with \cref{eq:hoffman-local}, this proves \cref{eq:positive-product-bound} near \(z\). Finitely many such neighborhoods cover \(Z\); on their compact complement, \(F\) has a positive minimum. Enlarging the constant completes the proof.
\end{proof}

Fix a lattice cube tiler \(U\Z^n\). Choose a compact matrix box \(\mathcal B_U\subset\mathrm{GL}_n^+(\R)\) centered at \(U\), with interior injective in lattice moduli. For
\[
 F_B=B[-1/2,1/2)^n,
\]
choose a fixed axis-parallel box \(R_U\) containing \(F_B\) for every \(B\in\mathcal B_U\). Define
\begin{equation}\label{eq:physical-multiplicity}
 N_B(x)=\sum_{k\in\Z^n}\one_{Q_n+Bk}(x),
 \qquad
 \Psi_U(B)=\int_{R_U}|N_B(x)-1|\,dx.
\end{equation}

\begin{proposition}[rectangular-atom formula]\label{prop:atom-formula}
The box \(\mathcal B_U\) has a finite polyhedral subdivision such that, on the closure \(P\) of each cell,
\begin{equation}\label{eq:atom-formula}
 \Psi_U(B)=
 \sum_{\alpha}c_\alpha
 \prod_{i=1}^n\ell_{\alpha i}(B),
 \qquad c_\alpha\in\Z_{>0},\quad \ell_{\alpha i}\ge0,
\end{equation}
where every \(\ell_{\alpha i}\) is an affine function on \(P\).
\end{proposition}

\begin{proof}
Only finitely many labels contribute on \(R_U\), uniformly in \(B\). Indeed, if \(x\in R_U\cap(Q_n+Bk)\), then \(k=B^{-1}(x-q)\) for some \(q\in Q_n\), and the compactness of \(\mathcal B_U\) bounds the right-hand side.

For each coordinate \(i\), order the finite affine endpoint family
\[
 (Bk)_i-\tfrac12,
 \qquad
 (Bk)_i+\tfrac12,
\]
together with the two endpoints of the \(i\)-th side of \(R_U\). The simultaneous weak endpoint orders define a finite polyhedral subdivision. On a cell closure, consecutive endpoint blocks cut that side of \(R_U\) into intervals of nonnegative affine lengths. Cartesian products of these intervals partition \(R_U\), up to null boundaries, into rectangular atoms. The set of cube labels containing the relative interior of an atom is constant on the cell. If that multiplicity is \(m_\alpha\), the atom contributes
\[
 |m_\alpha-1|\prod_{i=1}^n\ell_{\alpha i}(B).
\]
Atoms with coefficient zero are omitted. Coincident endpoints merely give a zero affine length, so the same formula holds on chamber boundaries.
\end{proof}

Let
\begin{equation}\label{eq:ambient-zero-set}
 \mathcal Z_U=
 \{B\in\mathcal B_U:Q_n+B\Z^n\text{ is a level-one lattice tiling}\}.
\end{equation}

\begin{proposition}[ambient fixed-window estimate]\label{prop:ambient-window}
After replacing \(\mathcal B_U\) by a smaller concentric box, there is \(C_U>0\) such that
\begin{equation}\label{eq:ambient-window-bound}
 \dist(B,\mathcal Z_U)^n\le C_U\Psi_U(B)
\end{equation}
for every \(B\) in the smaller box. Moreover,
\begin{equation}\label{eq:window-zero-set}
 \Psi_U(B)=0\quad\Longleftrightarrow\quad B\in\mathcal Z_U.
\end{equation}
\end{proposition}

\begin{proof}
Apply \cref{lem:positive-product} to \cref{eq:atom-formula} on every cell closure meeting the smaller box. If a closure meets its zero set, the lemma controls distance to that cellwise zero set; distance to the full zero set can only be smaller. If a closure does not meet the zero set, \(\Psi_U\) has a positive minimum on its compact portion in the smaller box, while distance to the nonempty set \(\mathcal Z_U\) is bounded. The finitely many cellwise constants therefore give \cref{eq:ambient-window-bound} once \cref{eq:window-zero-set} is established.

If \(B\in\mathcal Z_U\), then \(N_B=1\) almost everywhere and \(\Psi_U(B)=0\). Conversely, suppose \(\Psi_U(B)=0\). Since \(F_B\subset R_U\), one has \(N_B=1\) almost everywhere on \(F_B\). The unfolding identity for the \(B\Z^n\)-periodization gives
\[
 \int_{F_B}N_B(x)\,dx=\vol(Q_n)=1.
\]
The pointwise equality on \(F_B\) makes the same integral equal to \(\vol(F_B)=\det B\); hence \(\det B=1\). Periodicity propagates \(N_B=1\) from the fundamental domain to almost every point of \(\R^n\), proving that \(B\in\mathcal Z_U\).
\end{proof}

\begin{proposition}[two-sided comparison with overlap]\label{prop:window-overlap}
There are \(m_U,M_U>0\) such that
\begin{equation}\label{eq:window-overlap}
 m_U\Delta_n(B)\le \Psi_U(B)\le M_U\Delta_n(B)
\end{equation}
for every unimodular \(B\in\mathcal B_U\).
\end{proposition}

\begin{proof}
The fixed box \(R_U\) meets at most \(M_U\) translates of \(F_B\), uniformly on the compact matrix box. The function \(|N_B-1|\) is \(B\Z^n\)-periodic, so
\[
 \Psi_U(B)\le M_U\int_{F_B}|N_B-1|.
\]
For a nonnegative integer \(m\), \(|m-1|\le(m-1)^2\). The variance identity \cref{lem:variance} therefore gives
\[
 \int_{F_B}|N_B-1|
 \le\int_{F_B}(N_B-1)^2
 =\Delta_n(B),
\]
which proves the upper bound for \(\Psi_U\).

Compactness of \(\mathcal B_U\) also gives a uniform integer \(L_U\) such that \(N_B(x)\le L_U\) for every \(B\in\mathcal B_U\) and every \(x\in F_B\). Indeed, if \(x\in F_B\cap(Q_n+Bk)\), then \(k=B^{-1}(x-q)\) with \(q\in Q_n\), and the right-hand side ranges over a fixed compact set. Since \(F_B\subset R_U\) and \((m-1)^2\le L_U|m-1|\) for \(0\le m\le L_U\),
\[
 \Delta_n(B)
 =\int_{F_B}(N_B-1)^2
 \le L_U\int_{F_B}|N_B-1|
 \le L_U\Psi_U(B).
\]
Taking \(m_U=L_U^{-1}\) proves the other bound.
\end{proof}

\begin{theorem}[sharp local universal estimate]\label{thm:local-universal}
For every lattice cube tiler \(U\Z^n\), there are \(C_U,\varepsilon_U>0\) such that
\begin{equation}\label{eq:local-universal}
 \dist(B,\mathcal Z_U)
 \le C_U\Delta_n(B)^{1/n}
\end{equation}
for every \(B\in\mathrm{SL}_n(\R)\) with \(\|B-U\|<\varepsilon_U\). Consequently,
\begin{equation}\label{eq:kappa-bound}
 \kappa(U)\le n
\end{equation}
for every tiler.
\end{theorem}

\begin{proof}
Combine \cref{prop:ambient-window,prop:window-overlap}. The exact local exponent in \cref{thm:signed-valuations} is the largest order quotient along a Puiseux arc, so \cref{eq:local-universal} also gives \cref{eq:kappa-bound}.
\end{proof}

\begin{proof}[Proof of \cref{thm:main}]
The compact set \(\cT_n\) has a finite cover by injective basis charts in which \cref{thm:local-universal} holds. Matrix distance and the chosen quotient Riemannian distance are bi-Lipschitz on smaller chart closures. Properness of \(\cE_n\), together with \(\cE_n^{-1}(0)=\cT_n\), supplies \(\delta_n>0\) such that the sublevel \(\{\cE_n\le\delta_n\}\) lies in this finite cover. Taking the largest adjusted chart constant proves \cref{eq:intro-main}.

For sharpness, take the determinant-one cyclic family \(B_n(t,\ldots,t)\) of \cref{eq:cyclic-family}. The exact formula in \cref{thm:cyclic-overlap} gives
\[
 \cE_n(B_n(t,\ldots,t)\Z^n)
 =2|t|^n(1-|t|^{n-1}).
\]
Every total order has a backward edge on the directed Hamiltonian cycle, so the max-norm distance to the local branch union is at least \(|t|\). For sufficiently small \(t\), the injective chart contains every nearest lift of \(\cT_n\); local bi-Lipschitz equivalence therefore gives
\[
 \dist_{X_n}(B_n(t,\ldots,t)\Z^n,\cT_n)\asymp|t|.
\]
An estimate with exponent larger than \(1/n\) is impossible along this family.
\end{proof}

\begin{corollary}[exact global contact exponent]\label{cor:global-kappa}
Let \(\kappa_n=\max_{U\Z^n\in\cT_n}\kappa(U)\), with \(\kappa(U)\) as in \cref{eq:kappa-U}. Then
\[
 \kappa_n=n.
\]
\end{corollary}

\begin{proof}
The local universal estimate gives \(\kappa(U)\le n\) at every tiler. The cyclic arc at the standard lattice gives \(\kappa(I)\ge n\).
\end{proof}

\section{The sharp planar theorem}\label{sec:planar}

In dimension two the two global branches have representatives
\[
 U(a)=\begin{pmatrix}1&a\\0&1\end{pmatrix},
 \qquad
 L(b)=\begin{pmatrix}1&0\\b&1\end{pmatrix},
 \qquad a,b\in\R/\Z.
\]

\begin{proposition}[topology of the planar tiling locus]\label{prop:planar-topology}
The locus \(\cT_2\) is the union of two embedded circles meeting only at the standard lattice. In particular,
\[
 \cT_2\simeq S^1\vee S^1.
\]
\end{proposition}

\begin{proof}
\Cref{thm:nil-branches} identifies each branch with \(N_2(\R)/N_2(\Z)\simeq\R/\Z\). If \(U(a)\Z^2=L(b)\Z^2\), then
\[
 U(a)^{-1}L(b)=
 \begin{pmatrix}1-ab&-a\\b&1\end{pmatrix}
\]
is integral. Thus \(a,b\in\Z\), and both points are the standard lattice in their circle quotients.
\end{proof}

Fix \(x\in(0,1)\) and write a nearby determinant-one matrix as
\[
 B=\begin{pmatrix}p&q\\r&s\end{pmatrix},
 \qquad ps-qr=1,
 \qquad B\approx U(x).
\]
Only the three unoriented labels \(e_1,e_2,e_2-e_1\) contribute. In a neighborhood on which the remaining signs are fixed,
\begin{equation}\label{eq:planar-exact}
 \frac12\Delta_2(B)
 =(1-p)_+(1-|r|)
 +(1-q)(1-s)_+
 +(1-p+q)(1-s+r)_+.
\end{equation}
Put
\[
 a_1=1-p,\qquad a_2=1-s,\qquad c=r.
\]
The tangent determinant constraint at \(U(x)\) is
\begin{equation}\label{eq:planar-tangent}
 a_1+a_2+xc=0.
\end{equation}
The degree-one transverse initial form of half the energy is
\begin{equation}\label{eq:Phi}
 \Phi_x(a_1,a_2,c)
 =[a_1]_+ +(1-x)[a_2]_+ +x[a_2+c]_+.
\end{equation}

\begin{lemma}[planar contact norm]\label{lem:planar-norm}
On the plane \eqref{eq:planar-tangent}, the function \(\Phi_x\) vanishes only at the origin. Uniformly for \(x\in[\eta,1-\eta]\), it is comparable to the Euclidean norm, with constants depending only on \(\eta>0\).
\end{lemma}

\begin{proof}
Suppose \(\Phi_x=0\). Then \(a_1\le0\), \(a_2\le0\), and \(a_2+c\le0\). Write \(a_1=-u\) and \(a_2=-v\), with \(u,v\ge0\). The tangent constraint gives \(xc=u+v\), while \(a_2+c\le0\) gives \(c\le v\). Hence
\[
 u+v=xc\le xv<v
\]
unless \(v=0\). It follows that \(u=v=c=0\). The restriction of the continuous positively homogeneous function \(\Phi_x\) to the unit sphere in \eqref{eq:planar-tangent} has a positive minimum. Compactness in \(x\in[\eta,1-\eta]\) makes the minimum uniform.
\end{proof}

\begin{proposition}[linear rigidity on regular planar branches]\label{prop:planar-linear}
On every compact subarc of either branch avoiding the standard lattice,
\[
 \dist_{X_2}(\Lambda,\cT_2)\le C\,\cE_2(\Lambda).
\]
\end{proposition}

\begin{proof}
Formula \eqref{eq:planar-exact} has initial form \eqref{eq:Phi} on the determinant tangent plane. By \cref{lem:planar-norm}, that initial form is a norm on a complement of the branch tangent. The exact determinant equation differs from its tangent equation by a quadratic term. Compactness of the subarc and the fixed finite contact table therefore give a uniform linear lower bound for the energy in the transverse distance.
\end{proof}

\begin{theorem}[global planar rigidity]\label{thm:global-planar}
There are \(C,\delta>0\) such that
\begin{equation}\label{eq:global-planar}
 \dist_{X_2}(\Lambda,\cT_2)
 \le C\,\cE_2(\Lambda)^{1/2}
\end{equation}
whenever \(\cE_2(\Lambda)\le\delta\). The exponent \(1/2\) is optimal, and the bound improves to a linear one locally uniformly away from the standard lattice.
\end{theorem}

\begin{proof}
Properness implies that every sufficiently small energy sublevel lies in a prescribed neighborhood of the compact zero set. Near the standard lattice, \cref{thm:sharp-standard} with \(n=2\) gives \eqref{eq:global-planar}. Remove a smaller neighborhood of that crossing from the two branch circles. The remainder is compact and is covered by finitely many charts from \cref{prop:planar-linear}. Since \(0\le s\le1\) implies \(s\le s^{1/2}\), the estimates combine.

For optimality, take the determinant-one two-cycle family with both off-diagonal entries of size \(t\). Its distance to the branch crossing is comparable to \(|t|\), while \cref{thm:normal-form} gives energy comparable to \(t^2\).
\end{proof}

\section{Cyclic extremizers and spectral escape}\label{sec:spectral}

For \(q=(q_1,\ldots,q_n)\), put \(p(q)=\prod_iq_i\) and define
\begin{equation}\label{eq:cyclic-family}
 B_n(q)
 =I+\sum_{i=1}^{n-1}q_iE_{i,i+1}
 +q_nE_{n,1}+(-1)^np(q)E_{n,n}.
\end{equation}
The determinant of the uncorrected cyclic shear is \(1+(-1)^{n-1}p(q)\), while the cofactor of the last diagonal entry is one. Hence \(\det B_n(q)=1\).

\begin{theorem}[exact cyclic overlap]\label{thm:cyclic-overlap}
If \(\max_i|q_i|\le1/4\), then
\begin{equation}\label{eq:cyclic-overlap}
 \Delta_n(B_n(q))
 =2|p(q)|
 \left(1-\prod_{i=1}^{n-1}|q_i|\right).
\end{equation}
\end{theorem}

\begin{proof}
If some \(q_i\) vanishes, the shear support is acyclic and both sides vanish. Assume all \(q_i\ne0\), and put \(\delta=(-1)^np(q)\).

Only labels in \(\{-1,0,1\}^n\) contribute. Indeed, if \(m=\max_i|k_i|\ge2\), a row where that maximum is attained gives
\[
 |(B_n(q)k)_i|\ge m(1-1/4)>1,
\]
with the last row handled by the additional term of size at most \(4^{-n}\).

Let \(k\in\{-1,0,1\}^n\setminus\{0\}\) contribute and let \(j\) be its first nonzero index. Positivity of the overlap factors in rows \(i<n\) forces
\begin{equation}\label{eq:tail-recursion}
 k_{i+1}=-\operatorname{sgn}(q_i)k_i
 \qquad\text{whenever }k_i\ne0.
\end{equation}
Thus, up to overall sign, every possible contributor is a tail: it vanishes before \(j\), and its later entries are fixed by \eqref{eq:tail-recursion}.

For \(j\ge2\), the product of the first \(n-1\) overlap factors is
\[
 (1-|q_{j-1}|)\prod_{i=j}^{n-1}|q_i|.
\]
The last coordinate is \((1+\delta)k_n\). These tails contribute exactly when \(\delta<0\), with last factor \(|p(q)|\). For \(j=1\), the product of the first \(n-1\) factors is \(\prod_{i=1}^{n-1}|q_i|\), and the last factor is \(|q_n|-|p(q)|\) exactly when \(\delta>0\).

If \(\delta>0\), the one unoriented active label contributes
\[
 |p(q)|\left(1-\prod_{i=1}^{n-1}|q_i|\right).
\]
If \(\delta<0\), the unoriented tails contribute
\[
 |p(q)|\sum_{j=2}^n
 (1-|q_{j-1}|)\prod_{i=j}^{n-1}|q_i|
 =|p(q)|\left(1-\prod_{i=1}^{n-1}|q_i|\right),
\]
by telescoping. Opposite labels contribute equally, proving the formula.
\end{proof}

\begin{corollary}\label{cor:cycle-sharpness}
For \(q_1=\cdots=q_n=t\),
\[
 \Delta_n(B_n(t,\ldots,t))
 =2|t|^n(1-|t|^{n-1}),
\]
whereas \(\dist_\infty(B_n(t,\ldots,t),\cH_n)\ge|t|\).
\end{corollary}

\begin{proof}
The energy formula is \cref{thm:cyclic-overlap}. Every total order has at least one backward edge on the directed Hamiltonian cycle, so every Haj\'os branch deletes an entry of magnitude \(|t|\).
\end{proof}

Let \(\Lambda^*\) denote the Euclidean dual lattice and use the convention \(\sinc(0)=1\).

\begin{theorem}[spectral overlap identity]\label{thm:spectral-identity}
For every unimodular lattice,
\begin{equation}\label{eq:spectral-identity}
 \cE_n(\Lambda)
 =\sum_{\xi\in\Lambda^*\setminus\{0\}}
 \prod_{j=1}^n\sinc^2(\xi_j).
\end{equation}
Consequently, \(Q_n\) tiles by \(\Lambda\) if and only if every \(0\ne\xi\in\Lambda^*\) has a coordinate in \(\Z\setminus\{0\}\).
\end{theorem}

\begin{proof}
The Fourier coefficient of the periodization \(N_\Lambda\) at \(\xi\in\Lambda^*\) is
\[
 \widehat N_\Lambda(\xi)
 =\widehat{\one_{Q_n}}(\xi)
 =\prod_{j=1}^n\sinc(\xi_j).
\]
The constant coefficient is one. Parseval on \(\R^n/\Lambda\), followed by \cref{lem:variance}, proves \eqref{eq:spectral-identity}. Every summand is nonnegative, and \(\sinc(x)\) vanishes precisely at the nonzero integers. This proves the criterion, in agreement with the classical tiling--spectral framework \cite{Fuglede,IosevichPedersen,Kolountzakis2004}.
\end{proof}

For \(0\ne k\in\Z^n\), define the spectral summand
\begin{equation}\label{eq:Sk}
 S_k(q)=
 \prod_{j=1}^n
 \sinc^2\bigl((B_n(q)^{-T}k)_j\bigr).
\end{equation}

\begin{theorem}[fixed-mode law]\label{thm:fixed-mode}
Assume every \(q_i\ne0\) and \(\max_i|q_i|\to0\), with arbitrary relative rates. For each fixed \(0\ne k\in\Z^n\), exactly one of the following holds:
\begin{enumerate}[label=\textup{(\roman*)}]
\item if \(k\in\Z e_n\), then \(S_k(q)=0\) identically;
\item otherwise,
\begin{equation}\label{eq:fixed-mode}
 S_k(q)=|p(q)|^2(1+o(1)).
\end{equation}
\end{enumerate}
\end{theorem}

\begin{proof}
Write \(\xi=B_n(q)^{-T}k\). The system \(B_n(q)^T\xi=k\) is
\begin{align*}
 k_1&=\xi_1+q_n\xi_n,\\
 k_i&=\xi_i+q_{i-1}\xi_{i-1}\quad(2\le i<n),\\
 k_n&=(1+(-1)^np)\xi_n+q_{n-1}\xi_{n-1}.
\end{align*}
If \(k=me_n\), these equations give \(\xi_n=m\), so the last sinc factor vanishes.

Suppose \(k\notin\Z e_n\), and let \(S=\{i:k_i\ne0\}\). Follow the recurrence from one element of \(S\) to the next in cyclic order. If \(j\) precedes \(i\) in \(S\), and \(Q_{ji}\) is the product of the \(q\)-edges across that cyclic gap, then
\begin{equation}\label{eq:gap}
 \xi_i-k_i=\pm Q_{ji}k_j(1+o(1)).
\end{equation}
The determinant correction has full-cycle degree and is negligible against every proper-gap product. If \(S\) is a singleton, it supplies the full-cycle term; the only exceptional singleton is the already treated \(e_n\)-axis.

More explicitly, if \(k=me_j\) with \(j<n\), one circuit of the recurrence gives
\[
 \xi_j=m\bigl(1+(-1)^np(q)\bigr).
\]
Thus its unique active sinc factor has absolute value \(|p(q)|(1+o(1))\). For \(j=n\), the same circuit gives \(\xi_n=m\), which is precisely the exceptional case.

For \(m\in\Z\setminus\{0\}\) and \(\varepsilon\to0\),
\[
 |\sinc(m+\varepsilon)|
 =\frac{|\varepsilon|}{|m|}(1+o(1)).
\]
Coordinates outside \(S\) contribute \(1+o(1)\). Multiplication of \eqref{eq:gap} over \(i\in S\) telescopes the integer ratios \(|k_j/k_i|\), while the gap products use every \(q\)-edge exactly once. Thus the absolute Fourier coefficient is \(|p(q)|(1+o(1))\); squaring proves \eqref{eq:fixed-mode}.
\end{proof}

\begin{theorem}[spectral escape]\label{thm:spectral-escape}
Normalize the spectral summands to a probability measure
\[
 \mu_q(k)=
 \frac{S_k(q)}{\cE_n(B_n(q)\Z^n)},
 \qquad k\in\Z^n\setminus\{0\}.
\]
For every finite \(F\subset\Z^n\setminus\{0\}\),
\begin{equation}\label{eq:escape}
 \mu_q(F)\longrightarrow0.
\end{equation}
More precisely, if \(F^\circ=F\setminus\Z e_n\), then
\[
 \sum_{k\in F}S_k(q)
 =|F^\circ|\,|p(q)|^2(1+o(1)),
 \qquad
 \mu_q(F)=\frac{|F^\circ|}{2}|p(q)|(1+o(1)).
\]
\end{theorem}

\begin{proof}
Sum \cref{thm:fixed-mode} over the fixed finite set. By \cref{thm:cyclic-overlap}, the denominator is \(2|p(q)|(1+o(1))\).
\end{proof}

\begin{remark}
Every fixed mode is quadratic in the cycle monomial, but the complete positive spectral sum is linear in its absolute value. The sharp contact order is therefore carried by frequencies escaping to infinity; no finite-mode Hessian can detect it.
\end{remark}

\section{The poset normal form}\label{sec:scope}

Direct comparison of individual deck-contact summands with rooted monomials fails at a nonstandard tiler because cancellation may occur inside a complete affine boundary-row form. The fixed-window discrepancy removes this obstruction: its rectangular atoms are nonnegative, and its endpoint chambers split independently by matrix row. The following elementary rigidity principle converts that product structure and the already known zero set into the required analytic comparison.

Write \(\mathcal C(P)\) for the rooted \(P\)-circuits and set
\begin{equation}\label{eq:rooted-gauge-general}
 G_P(A)=
 \sum_i|a_{ii}|+
 \sum_{C\in\mathcal C(P)}\prod_{e\in C}|a_e|.
\end{equation}

\begin{theorem}[rooted contact comparison]\label{thm:rooted-contact}
For every cube-tiling lattice \(\Lambda=U\Z^n\), there are \(c_U,C_U,\varepsilon_U>0\) such that, in the coordinate \(B=(I+A)U\),
\begin{equation}\label{eq:rooted-contact-comparison}
 c_U G_{P_\Lambda}(A)
 \le\Delta_n(B)
 \le C_U G_{P_\Lambda}(A)
\end{equation}
whenever \(\det(I+A)=1\) and \(\|A\|<\varepsilon_U\).
\end{theorem}

\begin{lemma}[block-squarefree zero-set rigidity]\label{lem:block-squarefree}
Let \(V=V_1\oplus\cdots\oplus V_n\), and let \(F,G\) be nonnegative continuous functions on a neighborhood of the origin in \(V\). Suppose that a finite common subdivision near the origin consists of products \(K_1\times\cdots\times K_n\) of polyhedral cones. Suppose also that, on every product cone, each of \(F\) and \(G\) is comparable to a finite positive sum of block-squarefree products
\begin{equation}\label{eq:block-squarefree-sum}
 \sum_{\alpha}c_\alpha
 \prod_{i\in S_\alpha}\ell_{\alpha i}(x_i),
 \qquad c_\alpha>0,
\end{equation}
where \(S_\alpha\subset[n]\) and every \(\ell_{\alpha i}\) is linear and nonnegative on \(K_i\). If the germs \(F^{-1}(0)\) and \(G^{-1}(0)\) agree, then
\[
 F\asymp G
\]
near the origin.
\end{lemma}

\begin{proof}
Take a finite common refinement of the two product-conical subdivisions. For a row cone \(K_i\), let
\[
 L_i=K_i\cap(-K_i)
\]
be its lineality space. Every linear form nonnegative on \(K_i\) vanishes on \(L_i\): for \(y\in L_i\), nonnegativity at both \(y\) and \(-y\) forces \(\ell(y)=0\). Thus every model sum in \eqref{eq:block-squarefree-sum} factors through the pointed quotient \(K_i/L_i\). Choose a linear complement to \(L_i\), take a finite common pointed simplicial refinement of the resulting quotient cones, and pull it back to \(K_i\). It is therefore enough to work on one product of these refined cones.

On one such pointed simplicial quotient cone, write
\[
 \bar x_i=\sum_{v\in R_i}r_vu_v,
 \qquad r_v\ge0,
\]
where \(R_i\) is its ray set and \(\bar x_i\) is the class of \(x_i\) modulo \(L_i\). A linear form nonnegative on the original cone descends to the quotient and has an expansion
\[
 \ell_{\alpha i}(x_i)
 =\bar\ell_{\alpha i}(\bar x_i)
 =\sum_{v\in R_i}\bar\ell_{\alpha i}(u_v)r_v
\]
with nonnegative coefficients. Expanding \eqref{eq:block-squarefree-sum} therefore makes each function comparable to a positive squarefree monomial sum
\begin{equation}\label{eq:ray-hypergraph-sum}
 H_{\mathcal G}(r)=\sum_{E\in\mathcal G}d_E\prod_{v\in E}r_v,
 \qquad d_E>0,
\end{equation}
in which an edge \(E\) uses at most one ray from each row block.

If an expansion has a positive constant term, the corresponding function is bounded away from zero on the cone. Equality of zero germs then makes the same true for the other function, and comparability is immediate. Assume therefore that every edge is nonempty. Let \(T\) be any set of ray variables and choose a point, arbitrarily close to the origin, with \(r_v=0\) exactly for \(v\in T\) and all remaining ray coordinates positive. Positivity of the coefficients in \eqref{eq:ray-hypergraph-sum} gives
\[
 H_{\mathcal G}(r)=0
 \quad\Longleftrightarrow\quad
 T\cap E\ne\varnothing\ \text{for every }E\in\mathcal G.
\]
Thus the coordinate-zero set is determined exactly by those \(T\) that meet every hyperedge, namely the vertex covers of \(\mathcal G\). Comparability of each original function with its monomial model identifies their zero germs. Hence equality of the original zero germs implies that the two ray hypergraphs have exactly the same vertex covers.

The family of vertex covers determines its clutter of inclusion-minimal covers. Hypergraph blocker involution, applied once more, recovers from that clutter the original clutter of inclusion-minimal edges. The two hypergraphs consequently have the same minimal-edge clutter \(\mathcal M\).

Now restrict to \(0\le r_v\le1\). Every edge \(E\) contains some \(E_0\in\mathcal M\), and
\[
 \prod_{v\in E}r_v\le\prod_{v\in E_0}r_v.
\]
The positive terms indexed by \(\mathcal M\) give the reverse bound. Since all sums are finite and all their coefficients are positive, each monomial model is therefore comparable to
\[
 \sum_{E\in\mathcal M}\prod_{v\in E}r_v.
\]
Both functions are comparable to this same minimal-monomial sum on the refined cone. Pulling back across the lineality quotient proves comparability on the original product cone, and finiteness of the common subdivision makes the constants uniform near the origin.
\end{proof}

\begin{proof}[Proof of \cref{thm:rooted-contact}]
Put
\[
 W_U(A)=\Psi_U((I+A)U).
\]
Choose the matrix box in \cref{prop:atom-formula} to be a product of row boxes. Every endpoint used in its subdivision has the form
\[
 ((I+A)Uk)_i\pm\tfrac12
 =(Uk)_i+\sum_j a_{ij}(Uk)_j\pm\tfrac12,
\]
so an endpoint comparison in coordinate \(i\) depends only on the \(i\)-th row of \(A\). Comparisons whose endpoint difference is nonzero at \(A=0\) have a fixed sign after the neighborhood is reduced. Every remaining comparison wall passes through the origin and is given by a homogeneous linear equation in that row. The endpoint orders in distinct coordinates are chosen independently when the rectangular atoms are formed. Consequently the local subdivision is a common refinement of products
\[
 K_1\times\cdots\times K_n
\]
of rowwise polyhedral cones, exactly as required in \cref{lem:block-squarefree}.

On a cell closure, \cref{prop:atom-formula} writes \(W_U\) as a positive sum of rectangular-atom volumes. Each volume contains exactly one nonnegative affine interval gap from each coordinate row. A gap positive at the origin stays bounded above and below by positive constants, so it is a unit for comparability and may be absorbed into the coefficient. A gap vanishing at the origin becomes a nonnegative linear form on the corresponding row cone. After the unit factors are removed, every surviving monomial therefore contains at most one nonconstant factor from each row block. This is the block-squarefree structure required for \(W_U\).

Refine the same row cones by the signs of all coordinates \(a_{ij}\). On every refined cone, each absolute value in \eqref{eq:rooted-gauge-general} is a nonnegative linear form. The diagonal term \(|a_{ii}|\) belongs to row \(i\). A rooted circuit has distinct tails, and the tail of an edge \(i\to j\) is precisely the row index of \(a_{ij}\). Each rooted-circuit monomial therefore contains at most one coordinate from every matrix row. Hence \(G_{P_\Lambda}\) satisfies the other block-squarefree hypothesis.

It remains to identify the zero germs. In the fixed-window neighborhood, \cref{prop:ambient-window} gives
\[
 W_U(A)=0
 \quad\Longleftrightarrow\quad
 (I+A)U\Z^n\text{ is a level-one cube tiling}.
\]
The local affine Haj\'os atlas in \cref{prop:affine-atlas}, together with \eqref{eq:poset-arrangement}, identifies these nearby tilers exactly with \(A\in\cH(P_\Lambda)\). On the other hand, \cref{thm:rooted-ideal} identifies the zero germ of the rooted-circuit ideal, and hence of its positive absolute-value gauge:
\[
 G_{P_\Lambda}(A)=0
 \quad\Longleftrightarrow\quad
 A\in\cH(P_\Lambda).
\]
All equivalences are local at \(A=0\); after one common reduction of the neighborhood, the zero germs agree. The block-squarefree lemma therefore yields
\begin{equation}\label{eq:window-rooted-comparison}
 W_U(A)\asymp G_{P_\Lambda}(A)
\end{equation}
in the full ambient matrix neighborhood.

Now restrict this comparison to the smooth slice \(\det(I+A)=1\). On that slice, \cref{prop:window-overlap} replaces the fixed-window defect by overlap:
\[
 W_U(A)\asymp\Delta_n((I+A)U).
\]
Combining this with \eqref{eq:window-rooted-comparison} proves \eqref{eq:rooted-contact-comparison}.
\end{proof}

\begin{corollary}[stratum contact exponent]\label{cor:poset-contact-exponent}
For every cube-tiling lattice \(\Lambda=U\Z^n\),
\[
 \kappa(U)=c(P_\Lambda),
\]
and the optimal local repair exponent is \(1/c(P_\Lambda)\).
\end{corollary}

\begin{proof}
The comparison theorem and \cref{cor:algebraic-repair} give the exponent \(1/c(P_\Lambda)\). It remains to prove sharpness on the determinant-one slice. Choose a rooted circuit \(C\) of maximum size \(m=c(P_\Lambda)\), put the same parameter \(t\) on its variable edges, and set all other normal entries to zero. Inclusion-minimality of \(C\) implies that the rooted gauge has order \(m\), while \cref{thm:poset-bottleneck} gives distance of order one in \(t\), that is,
\[
 G_{P_\Lambda}(A(t))\asymp |t|^m,
 \qquad
 \dist(A(t),\cH(P_\Lambda))\asymp |t|.
\]

Let \(D(t)=\det(I+A(t))\). Every nonidentity determinant monomial is supported on the off-diagonal cycles of a nontrivial permutation. If such a monomial of degree strictly less than \(m\) were supported on the chosen circuit edges, one of its permutation cycles would use a nonempty proper subset of \(C\). 
That subset is already a forbidden variable support: it is a directed cycle even without the fixed $P$-edges. Therefore it contains a rooted $P_\Lambda$-circuit.
This contradicts the inclusion-minimality of \(C\). Hence
\[
 D(t)=1+O(t^m).
\]

Multiply one row of \(I+A(t)\) by \(D(t)^{-1}\). The resulting analytic matrix arc has determinant one. Since \(D(t)^{-1}=1+O(t^m)\), the normalization changes the selected off-diagonal entries only in higher order and introduces a diagonal defect of order at least \(m\). It follows that the normalized arc still has branch distance comparable to \(|t|\) and rooted gauge comparable to \(|t|^m\). By \cref{thm:rooted-contact}, its overlap is also comparable to \(|t|^m\). The arc characterization in \cref{thm:signed-valuations} proves sharpness and gives
\[
 \kappa(U)=m=c(P_\Lambda).
\]
\end{proof}

\section{Conclusion}

The sharp quantitative Haj\'os theorem gives a uniform stability law for the complete lattice cube-tiling locus: overlap energy controls distance with exponent \(1/n\), and cyclic degenerations show that this global exponent cannot be improved. The fixed-window multiplicity defect supplies the positive geometric observable behind the upper bound. Its decomposition into rectangular atoms reduces the all-strata estimate to products of at most \(n\) affine endpoint gaps.

The local geometry is finer. Incidence posets organize the branches through each tiler, and rooted circuits encode their transverse contacts. Rowwise endpoint subdivision and block-squarefree zero-set rigidity bypass the signed cancellations present in individual deck-contact terms, giving a two-sided comparison between overlap and the rooted-circuit gauge at every tiler. Consequently the local contact order is exactly \(c(P_\Lambda)\), the length of the longest rooted circuit, and the optimal local repair exponent is \(1/c(P_\Lambda)\).

Thus the global worst-case exponent \(1/n\) and the stratum-dependent exponent \(1/c(P_\Lambda)\) describe complementary scales of the same stability theory. 
Together, the overlap observable, the poset stratification, and the rooted-circuit normal form provide a quantitative bridge among (i) lattice tilings; (ii) semialgebraic stability; (iii) combinatorial incidence structure; and (iv) local contact geometry.

\section*{Data Availability Statement}

The verifier is included in the Supplementary Materials accompanying this record.

The accompanying Math Verification Suite records the dependency graph of the results and provides independent finite consistency checks for the exact computations used in the proof. It is supplementary to, and not a substitute for, the analytic arguments in the article.

\end{document}